\documentclass[12pt,a4paper]{article}

\usepackage{amssymb}

\newcommand{\vs}{\vspace}
\newcommand{\hs}{\hspace}
\newcommand{\NN}{\nonumber}
\newcommand{\NI}{\noindent}
\newcommand{\F}{\frac}      
\newcommand{\T}[1]{{\tilde #1}}
\newcommand{\B}[1]{{\bar #1}}
\newcommand{\C}[1]{{\cal #1}}
\renewcommand{\P}{\partial}
\newcommand{\DS}{\displaystyle}

\renewcommand{\SS}{\scriptstyle}
\newcommand{\SSS}{\scriptscriptstyle}
\newcommand{\G}{\left}
\newcommand{\D}{\right}
\newcommand{\GA}{\langle}
\newcommand{\DA}{\rangle}
\newcommand{\Coben}{\mathbb{C}^{(p+q)\cdot p~*}}
\newcommand{\Cinv}{C_0^\infty (\mathbb{C}^{(p+q)\cdot p~*})}
\newcommand{\Cinvl}{C_0^\infty (\mathbb{C}^{(p+q)\cdot p~*})[[\lambda]]}
\newcommand{\Chom}{C_{hom}^\infty (\mathbb{C}^{(p+q)\cdot p~*})}
\newcommand{\Choml}{C_{hom}^\infty (\mathbb{C}^{(p+q)\cdot p~*})[[\lambda]]}
\newcommand{\stt}{{\tilde{*}}}
\newcommand{\Ad}{\dagger}
\renewcommand{\phi}{\varphi}
\newcommand{\tr}{\mbox{\rm tr}}
\newcommand{\otr}{{\otimes r}}
\newtheorem{lemma}{Lemma}[section]
\newtheorem{theorem}{Proposition}[section]
\newtheorem{corollary}{Corollary}[section]
\newtheorem{defin}{Definition}[section]
\newenvironment{definition}{\begin{defin} \rm}{\end{defin}} 

\begin{document}

\title{A Star Product for Complex Grassmann Manifolds}

\author{\bf Joachim Schirmer\thanks
{e-mail: Joachim.Schirmer@physik.uni-freiburg.de}}

\date{September 12, 1997}

\maketitle

\centerline{FR-THEP 21}

\begin{abstract}

We explicitly construct a star product for the complex Grassmann manifolds 
using the method of phase space reduction. Functions on 
$\mathbb{C}^{(p+q)\cdot p~*}$, the space of $(p+q)\times p$ matrices of 
rank p, invariant under the right action of $Gl(p,\mathbb{C})$ can be 
regarded as functions on the Grassmann manifold $G_{p,q}(\mathbb{C})$, but 
do not form a subalgebra whereas functions only invariant under the unitary  
subgroup $U(p)\subset Gl(p,\mathbb{C})$ do. The idea is to construct a 
projection from $U(p)$- onto $Gl(p,\mathbb{C})$-invariant functions, whose 
kernel is an ideal. This projection can be used to define a star-algebra 
on $G_{p,q}(\mathbb{C})$ onto which this projection acts as an 
algebra-epimorphism.

\end{abstract}

\newpage

\section{Introduction}
The idea of deformation quantization is to modify the multiplication on the 
algebra of smooth complex valued functions in a noncommutative way, such 
that in zeroth order the normal multiplication is preserved, the commutator 
equals in first order the Poisson bracket, the constant function with value 
one acts again as one in the new algebra and the pointwise complex 
conjugation remains an antilinear involution. This concept of a 
quantization, first introduced by Berezin \cite{Ber75}, later as a program 
formulated by F. Bayen, M. Flato, c. Fronsdal, A. Lichnerowicz D. 
Sternheimer\cite{BFFLS78}, is obviously founded on the existence of this 
kind of a deformed product. This nontrivial question was settled for 
symplectic manifolds by M. DeWilde and P.B.A. Lecomte \cite{DWL83}, 
independently by B. Fedosov \cite{Fed85, Fed94} and by H.  Omori, Y. Maeda 
and A. Yoshioka \cite{OMY91}. 

Explicit examples of star products for nontrivial phase spaces of physical 
interest are still rare. Although Grassmann manifolds $G_{p,q}(\mathbb{C})$  
may be seen as the reduced phase space of a system of $(p+q)\cdot p$ 
harmonic 
oscillators where certain energy and angular momentum sums are fixed, it is  
not really a physical example. So the aim of this paper is principally to 
show that phase space reduction may be used to construct star 
products on more difficult symplectic manifolds. The methods used might be 
relevant as well to the possible deformation quantization of constrained 
systems.

C. Moreno \cite{Mor86} has already stated a recursion formula for a star 
product on nonexceptional K\"{a}hler symmetric spaces. M. Cahen, S. Gutt and
 J. 
Rawnsley \cite{CGR93} defined a star product on every compact Hermitian 
symmetric space and proved the convergence for the algebra of 
representative functions. Regarding $G_{p,q}(\mathbb{C})$ as a coadjoint 
orbit of  $U(p+q)$ Karabegov \cite{Kar94} has given a method of deformation 
quantization  in terms of representative functions which probably could be 
specialized explicitly in this case. The existence of star products on 
coadjoint orbits was treated as well by D. Arnal, J. Ludwig and M. Masmoudi 
\cite{ALM94}. For the special case of the projective space a closed formula 
for a star product of Wick type was given by M. Bordemann {\it et al.} 
\cite{BBEW95}. It was constructed as well by phase space reduction and will 
be generalized in the present paper. Whereas there an equivalence 
transformation of the Wick product on $\mathbb{C}^{q+1}\backslash\{0\}$ 
onto a new star product was constructed which could be projected to 
$\mathbb{C}P^q$ trivially, in the case of Grassmann manifolds 
$G_{p,q}(\mathbb{C})$ we will  directly head for the projection, because a 
similar equivalence transformation does not exist since the group 
$Gl(p,\mathbb{C})$ is noncommutative for $p>1$. 

The paper will be organized as follows: In the next section we will show 
how to construct a projection of functions on $\mathbb{C}^{(p+q)\cdot p~*}$  
onto functions on $G_{p,q}(\mathbb{C})$ whose kernel is an ideal with 
respect to the Wick star product on $C^\infty(\mathbb{C}^{(p+q)\cdot 
p~*})$. Then we will calculate this projection for the well known case of 
the projective space $\mathbb{C}P^q=G_{1,q}(\mathbb{C})$, thereby verifying 
the result of Bordemann {\it et. al} \cite{BBEW95} in a very simple way. We 
use this as a guideline in the fourth section where the formula for the 
star product on Grassmann manifolds is proved. The arising coefficients are 
treated in the last section, where we will state some open questions as 
well. Throughout this paper we use Einstein's summation convention.  

\section{The decomposition of the space  of invariant \\ 
         functions}\label{secdec} Our aim is to construct a star product on 
Grassmann manifolds $G_{p,q}(\mathbb{C})$ by making use of the Wick product 
on $\mathbb{C}^{(p+q)\cdot p~*}$, the space of complex $(p+q)\times p$ 
matrices of rank $p$, from which the Grassmann manifolds can be obtained by 
projection $\pi$. Unfortunately the homogeneous functions on 
$\mathbb{C}^{(p+q)\cdot p~*}$ -- functions invariant under the 
$Gl(p,\mathbb{C})$-action, which are the bijective image of functions on 
$G_{p,q}(\mathbb{C})$ by the pullback $\pi^*$ -- do not form a subalgebra. 
Yet the set of invariant functions -- we reserve this term for functions 
invariant under the unitary subgroup $U(p)$ of $Gl(p,\mathbb{C})$ -- is a 
subalgebra which we will denote by $(C_0^\infty (\mathbb{C}^{(p+q)\cdot 
p~*}) [[\lambda]],\tilde{*})$. So the Wick star product of two homogeneous 
functions will only be an invariant function and we are forced to search 
for a certain "projection" $P:=C_0^\infty (\mathbb{C}^{(p+q)\cdot p~*}) 
[[\lambda]]\longrightarrow C_{hom}^\infty (\mathbb{C}^{(p+q)\cdot p~*}) 
[[\lambda]]$, which satisfies $P(P(f\:\stt\:g)\:\stt\:h)= P(f\:\stt\: 
P(g\:\stt\:h))$  for any three homogeneous functions $f,g,h$. Then we can 
obviously equip $C_{hom}^\infty (\mathbb{C}^{(p+q)\cdot p~*}) [[\lambda]]$ 
with an algebra structure by defining $f*g:=P(f\:\stt\:g)$. For the unity 
function we would like the equation $f=f*1=P(f\:\stt\:1)=Pf$ to hold, thus 
$P\!\mid_{ C_{hom}^\infty (\mathbb{C}^{(p+q)\cdot p~*}) [[\lambda]]}=  
id\!\mid_{ C_{hom}^\infty (\mathbb{C}^{(p+q)\cdot p~*}) [[\lambda]]} $ and 
$ P^2=P$, so the name projection for $P$ is justified. We can split the 
algebra $C_0^\infty (\mathbb{C}^{(p+q)\cdot p~*}) [[\lambda]]=PC_0^\infty 
(\mathbb{C}^{(p+q)\cdot p~*}) [[\lambda]]\oplus \mbox{ker} P$ and notice 
that if $P$ satisfies the above rules and is thus an algebra epimorphism of 
$(C_0^\infty (\mathbb{C}^{(p+q)\cdot p~*}) [[\lambda]], 
\stt)$ onto $(C_{hom}^\infty (\mathbb{C}^{(p+q)\cdot p~*}) 
[[\lambda]],*)$, obviously ker$P$ must be a $\stt$-ideal $\C{J}_{\stt}$. 
This condition is in fact an equivalent formulation, the projection $P$ is 
an algebra epimorphism if and only if ker$P$ forms a $\stt$-ideal 
$\C{J}_{\stt}$. The decomposition of the space of invariant functions 
shall now be constructed. 

The complex Grassmann manifold $G_{p,q}(\mathbb{C})$  is the space of 
equivalence classes of the complex $(p+q)\times q$ matrices of rank $p$ 
where 
the equivalence class is defined as follows: 
\[ \NN  z_1 \sim z_2 \qquad  
\Longleftrightarrow \qquad \exists g\in Gl(p,\mathbb{C})\quad z_2=z_1g 
\qquad 
z_1,z_2\in \mathbb{C}^{(p+q)\cdot p~*} \]
\begin{equation}\label{pi}
\pi:\mathbb{C}^{(p+q)\cdot p~*}\longrightarrow G_{p,q}(\mathbb{C}),\quad 
z\longmapsto z/\!\!\sim  
\end{equation}
The Grassmann manifold $G_{1,q}$ is the projective space $\mathbb{C}P^q$. 
Grassmann manifolds can also be regarded as the reduced phase space 
associated to the right action $\Phi$ of  the unitary group $U(p)$:
\[ \Phi:\Coben \times U(p) \longrightarrow \Coben ,\quad (z,U)
\longmapsto \Phi_U z =zU  \] 
This action is symplectic with respect to the K\"ahler symplectic form 
\begin{equation}\label{omega}      
  \omega = \F{i}{2} \mbox{tr}\, dz\wedge dz^\Ad = \F{i}{2} 
  \P\B{\P}\,\mbox{tr}\, z^\Ad z = \F{i}{2}\P\B{\P}\,\mbox{tr}\, x \qquad 
  x:=z^\Ad z 
\end{equation}
and there is an Ad$^*$-equivariant momentum mapping
\begin{equation}\label{J}
  J:\mathbb{C}^{(p+q)\cdot 
  p~*}\longrightarrow u(p)^*,\qquad J=\F{i}{2}z^\Ad z=\F{i}{2} x 
\end{equation}
where $u(p)^* $ denotes the dual of the space of antihermitian matrices. 
Since for any positive definite hermitian $\mu\in H^+(p)$ there exists a 
representative of $z/\!\!\sim$ on $ J^{-1}(\F{i}{2}\mu)$, namely $z(z^\Ad 
z)^{-\F{1}{2}}\mu^{\F{1}{2}}$, the Grassmann manifold $G_{p,q}$ is the 
quotient $J^{-1}(\F{i}{2}\mu)/\mu^{-\F{1}{2}}U\mu^{\F{1}{2}}$, and this is 
the reduced phase space if and only if $\mu$ is a positive scalar multiple 
of the unit matrix, which is always the case for $p=1$. So we choose in the   
following $\mu\in\mathbb{R}^+$ and elucidate the situation by the usual 
phase space reduction diagramm \cite{AM85}{\vspace{2em}}:

\unitlength 0.80mm
\linethickness{0.4pt}
\begin{picture}(85.00,39.00)
(00.00,57.00)
\put(80.00,100.00){\makebox(0,0)[cc]{$\mathbb{C}^{(p+q)\cdot p~*}$}}
\put(80.00,60.00){\makebox(0,0)[cc]{$G_{p,q}(\mathbb{C})$}}
\put(80.00,94.00){\vector(0,-1){28.00}}
\put(85.00,80.00){\makebox(0,0)[cc]{$\pi$}}
\put(55.00,80.00){\makebox(0,0)[cc]{$J^{-1}(\F{i}{2}\mu)$}}
\put(60.00,76.00){\vector(3,-2){15.00}}
\put(58.00,84.00){\vector(3,2){15.00}}
\put(77.00,94.00){\vector(-3,-2){15.00}}
\put(64.00,91.00){\makebox(0,0)[cc]{$i$}}
\put(71.00,87.00){\makebox(0,0)[cc]{$\B{\rho}$}}
\put(65.00,70.00){\makebox(0,0)[cc]{$\pi_\mu$}}
\end{picture}

\NI Additionally to the general case there exist the 
direct projections $\pi$ and $\B{\rho}$. Identifying 
$\mathbb{C}^{(p+q)\cdot p~*}$ with $J^{-1}(\F{i}{2}\mu)\times H^+$ by the 
global diffeomorphism 
\begin{equation}\label{z-xzeta}
  z\longmapsto(\zeta,x)=(z(z^\Ad z)^{-\F{1}{2}}\mu^{\F{1}{2}},z^\Ad z)  
\end{equation}
the mappings $i$ and $\B{\rho}$ adopt the explicit form
\begin{equation}\label{}
  i(\zeta)=(\zeta,\mu)\qquad \B{\rho}(\zeta,x)=\B{\rho}(z)=\zeta=z(z^\Ad 
  z)^{-\F{1}{2}}\mu^{\F{1}{2}}~, 
\end{equation}
and an invariant function $f\in \Cinv $ satisfies
\begin{equation}\label{finvariant}
  f(z)=f(\zeta,x)=f(zU)=f(\zeta U,U^\Ad x U)~.
\end{equation}
Looking for an appropriate choice of the star ideal and the associated 
projection we consider the trivial case of the pointwise product as a 
guideline. There is obviously a 1:1 mapping between invariant functions on 
the orbit $J^{-1}(\F{i}{2}\mu)$ and functions on $G_{p,q}(\mathbb{C})$, i.e. 
$\pi_{\mu}^*$ is invertible for invariant functions. All invariant 
functions on the orbit may be obtained by restriction of invariant 
functions on $\mathbb{C}^{(p+q)\cdot p~*}$ and thus there is a simple 
projection $P_\cdot$ 
\begin{equation}\label{Pmal}
  P_{\cdot}:\Cinv\longrightarrow\Chom,\quad 
  P_{\cdot}:=\pi^*\circ{\pi_\mu^*}^{-1}\circ i^*=\B{\rho}^*\circ i^*\quad.  
\end{equation}
We can easily decompose an invariant function into its homogeneous and its 
ideal part, thereby describing the ideal:
\begin{eqnarray}
\NN   \big[(1-P_\cdot)F\big](z) & = & F(\zeta,x)-F(\zeta,\mu) =
\G\GA\int_0^1dt \F{\P F}{\P x}\G(\zeta,tx 
                      +   (1-t)\mu \D), x-\mu\D\DA 
\end{eqnarray}
So ker$P_\cdot = (1-P_\cdot)\Cinv$ is generated by $x-\mu$:
\[ \mbox{ker}P=\C{J_{\cdot}}=\GA x-\mu \DA =\{F\mid  F=\GA x-\mu, G\DA,\quad
G\in C_{eq}^{\infty}\G(\mathbb{C}^{(p+q)\cdot p~*},u(p)\D)\}  \]            
Note  that in our case of invariant functions we can restrict the domain of  
functions paired with $x-\mu$ onto Ad-equivariant functions since for 
\begin{equation}\label{H}
  H:\Cinv\longrightarrow C_{eq}^{\infty}(\mathbb{C}^{(p+q)\cdot p~*}), 
  \quad H(F)(\zeta,x):=\int_0^1dt \F{\P F}{\P x}(\zeta,tx+(1-t)\mu) 
\end{equation}
holds 
\begin{equation}\label{Hequiv}
  H(F\circ\Phi_U)=\mbox{Ad}_{U^{-1}}\circ H(F) \circ \Phi_U = H(F)~.
\end{equation}
The decomposition 
\begin{equation}\label{Fdecomp}
  F=P_\cdot F+(1-P_\cdot)F = \B{\rho}^* i^*F + \GA H(F), x-\mu\DA
\end{equation}
is unique whereas the function $H$ is not, unless $p=1$, because we can add 
to $H(F)$ any function orthogonal to $x-\mu$. We have chosen $H$ here to be  
the integral along  the straight line between $(\zeta,x)$ and 
$(\zeta,\mu)$, but any other curve connecting these points will do as well.  
For the solution a convenient choice of $H(F)$ will be important as well as 
the proof that the decomposition is unique, such that the result is 
independent of the choice of $H(F)$. 

Analogously we would like to define a star ideal $\C{J}_{\stt}$ generated 
by $x-\mu$. First of all we have of course to extend our considerations to 
power series in a formal parameter $\lambda$ with coefficients in the 
respective function spaces. The ideal 
\[ 
\C{J}_{\stt} :=\GA x-\mu\DA_\stt = \{F\mid = \GA G\stackrel{\stt}{,}x-\mu
\DA = G^{i}_{~j}\:\stt\: (x-\mu)^{j}_{~i},\quad G\in C_{eq}^{\infty} 
(\mathbb{C}^{(p+q)\cdot p~*},u(p))[[\lambda]]\} 
\] 
may be defined independently 
of a global basis for $\mathbb{C}^{(p+q)\cdot p~*}$ and is indeed a 
twosided ideal because the Wick star product is strongly $U(p)$-invariant 
\begin{equation}\label{Ibeidseitig}
  f\:\stt\: J -J\: \stt\: f = \F{i\lambda}{2}\{f,J\} \qquad\forall f\in 
  C^\infty(\mathbb{C}^{(p+q)\cdot p~*})~, 
\end{equation}
while $U(p)$ is compact -- note that $G$ is only equivariant. Now we can 
formulate the important
\begin{lemma}\label{Zerlegung}
There exists a unique decomposition of invariant functions
\[ \Cinvl = \Choml \oplus \C{J}_{\stt} \] 
into a homogeneous and a $\stt$-ideal part, where the ideal is spanned by 
$x-\mu$. 
\end{lemma}
As an immediate consequence of the proposition and the remarks in the 
introduction we obtain the following
\begin{corollary}
There exists a star product on $G_{p,q}(\mathbb{C})$ and an algebra 
epimorphism 
\[ {\pi^{*}}^{-1}\circ P: (\Cinvl,\stt )\longrightarrow (C^{\infty}  
\G(G_{p,q}(\mathbb{C})\D)[[\lambda]],*)~. \] \vspace{-8.7ex}
\end{corollary}{\hspace*{\fill} $\square$}

\vspace{1ex}
\NI{\bf Proof:} We start with the point multiplicative decomposition and 
start 
deforming the usual product in the following way:
\begin{eqnarray}
 F& = & \B{\rho}^*i^*F + \GA H(F), x-\mu\DA  ~~=~~  \B{\rho}^*i^* F+\GA 
 H(F)\stackrel{\stt}{,} x-\mu\DA +\lambda \Delta 
  (F) \\ \label{Delta} \Delta(F) &:=& \F{1}{\lambda}\G (\GA 
  (H(F)\stackrel{\stt}{,} x-\mu \DA - \GA      H(F),x-\mu \DA \D) 
\end{eqnarray}
$\Delta(F)$ is a power series in $\lambda$ since the coefficient of 
$\lambda^{-1}$ vanishes and we can proceed decomposing $\Delta(F)$. Using 
induction, keeping in mind the linearity of $\B{\rho}^*,~ i^*,~H$ and 
$\Delta$ 
yields:
\begin{equation}\label{F*decomposed}
  F=\B{\rho}^*i^*\G(\F{1}{1+\lambda\Delta } F\D)+\G\GA H\G(\F{1} 
  {1+\lambda\Delta } F\D)  \stackrel{\stt}{,} x-\mu\D\DA 
\end{equation}
So there is a decomposition, but since $H$ is not uniquely defined, 
$\Delta$ is not either, so uniqueness might be lost during the deformation 
process. To insure uniqueness assume a function $F$ in the $\stt$-ideal 
$\C{J}_\stt$ is simultaneously homogeneous, i.e. 
\[ \B{\rho}^* \phi=F=\GA G\stackrel{\stt}{,} x-\mu \DA~,\quad\qquad G\in 
C^\infty_{eq} 
\G(\mathbb{C}^{(p+q)\cdot p~*},u(p)\D)[[\lambda]] \quad \phi \in 
C_0^\infty\G(J^{-1}(\F{i}{2}\mu)\D)[[\lambda]]
\] 
It has to be shown that $F$ vanishes identically, thus both sides of the 
above equation. We can make use of the fact, that the series for the Wick  
product with $x$ breaks off after the first term: 
\begin{equation}
\NN    \GA G\stackrel{\stt}{,} x\DA = \GA G, x\DA+\lambda \T{M}_1\GA G, x\DA, 
       \quad \T{M}_1\GA G,x\DA := z^A_{~i}\F{\P}{\P z^A_{~j}}G^i_{~j}\quad 
       G\in C^\infty_{eq}\G(\mathbb{C}^{(p+q)\cdot p~*},u(p)\D) 
\end{equation}                                  
Let $\varphi_{\nu }\in C^\infty_0(J^{-1}(\F{i}{2}\mu))$ and $G_{(\nu )}\in 
C^\infty_{eq}\G(\mathbb{C}^{(p+q)\cdot p~*},u(p)\D)$ be the coefficients of 
$\varphi$ and $G$ respectively:
\[ \B{\rho}^* \phi = \sum_{\nu =0}^{\infty} \lambda^\nu  \B{\rho}^* \phi_\nu 
 = \sum_{\nu =0}^\infty \lambda^\nu  \bigg (\GA G_{(\nu 
 )},x-\mu\DA+\T{M}_1\GA G_{(\nu - 1)}, x\DA \bigg) \qquad G_{(-1)} 
\equiv 0\]    
Considering all orders separately we obtain after restriction to  $x=\mu$ 
\begin{equation}\label{phires}
  \phi_\nu =i^* \T{M}_1\GA G_{(\nu  -1)},x-\mu\DA\qquad \phi_0=0 ~.
\end{equation}
The strategy is the following: Assume $\phi_\nu  =0 $ for $\nu  < n$. 
Consequently it holds 
\begin{equation}\label{M1M0}
  \T{M}_1\GA G_{(\nu  -1)},x-\mu\DA = -\GA G_{(\nu )},x-\mu\DA  \qquad 
  \mbox{for } \nu  < n~.
\end{equation}
Express $\phi_n=i^*\T{M}_1\GA G_{(n -1)},x\DA $ by derivatives of $\GA 
G_{(n 
-1)},x-\mu\DA $ and use the $n -1$ equations (\ref{M1M0}) successively 
in order to write $\phi_n$ as derivatives of $\GA G_{(0)},x-\mu\DA$ which 
vanishes identically. Then $\phi_n=0$ and by induction $\phi_\nu 
=0$ for all $\nu $ what completes the proof. 

Obviously we have to consider derivatives at $J^{-1}(\F{i}{2}\mu)$ 
transversally to this  submanifold. Thus in what follows we set $z=\zeta y$ 
and we investigate derivatives with respect to $y$ at the point $y=1$. In a 
tedious, but simple calculation one verifies the following formulas.  
\begin{eqnarray}
  \F{\P^r}{\P y^{i_1}_{~j_1}\ldots\P y^{i_r}_{~j_r}}\Bigg|_{y=1} \bigg(\GA 
  G_{(\nu )},x-\mu\DA(\zeta y) \bigg) & = & \mu\sum_{s=1}^r \F{\P^{r-1} 
  {G_{(\nu )}}^{j_s}_{~i_s}}{\P z^{A_1}_{~j_1}{\stackrel{\SS 
  s}{\SSS\wedge}\atop\DS\cdots} \P z^{A_r}_{~j_r}}(\zeta) 
  \;\zeta^{A_1}_{~i_1} {\stackrel{\SS 
  s}{\SSS\wedge}\atop\DS\cdots}\zeta^{A_r}_{~i_r}\\ \NN \F{\P^r}{\P   
  y^{i_1}_{~j_1}\ldots\P y^{i_r}_{~j_r}}\Bigg|_{y=1} \bigg(\T{M}_1\GA 
  G_{(\nu )},x-\mu\DA(\zeta y) \bigg) & = &   
  \delta^{i_{r+1}}_{~j_{r+2}}\delta^{i_{r+2}}_{~j_{r+1}}\F{\P^{r+1} 
  {G_{(\nu )}}^{j_{r+2}}_{~i_{r+2}}}{\P z^{A_1}_{~j_1}\ldots\P 
  z^{A_{r+1}}_{~j_{r+1}}}(\zeta)\; 
  \zeta^{A_1}_{~i_1}\ldots\zeta^{A_{r+1}}_{~i_{r+1}}\\ && 
  \hs{-10ex}+\sum_{s=1}^r \delta^{i_{r+1}}_{~j_{r+1}}\F{\P^r {G_{(\nu 
  )}}^{j_s}_{~i_{r+1}}}{\P z^{A_1}_{~j_1}\ldots \P z^{A_s}_{~j_{r+1}}\ldots 
  \P z^{A_r}_{~j_s}}(\zeta)\; \zeta^{A_1}_{~i_1}\ldots \zeta^{A_r}_{~i_r} 
\end{eqnarray}

Before proceeding with the proof let us summarize some notation for 
conjugacy classes of the symmetric group. By the $r$-tuple 
$\alpha=(\alpha_1,\ldots,\alpha_r)$ we denote the conjugacy class in $S_r$ 
having $\alpha_i~i$-cycles in the cycle representation, by $Con(S_r)$ the 
set of all conjugacy classes  Obviously it holds $\sum i\alpha_i=r$ and we 
define $|\alpha|=\alpha_1+\ldots+\alpha_r$ as usual. The number of elements 
in the  conjugacy class $\alpha$ is 
\begin{equation}\label{h_alpha}
  h_\alpha = \F{r!}{1^{\alpha_1} \alpha_1 !\cdot\ldots\cdot 
  r^{\alpha_r}\alpha_r !}~. 
\end{equation}
Since $i-1$ transpositions are necessary to compose an $i$-cycle from 
1-cycles (resp. to decompose it into 1-cycles) the lowest number of 
transpositions, denoted by \#, that is needed to compose a certain 
permutation $\sigma\in \alpha$ from the identity (or vice versa to 
decompose it into the identity) is only dependent on the conjugacy class 
$\alpha$ and it holds 
\begin{equation}\label{step}
  \# \alpha=\#[ \sigma]=  \sum_{i=1}^r(i-1)\alpha_i=r-|\alpha|
\end{equation}
We abbreviate by $D:=\G(\F{\P}{\P y^{i}_{~j}}\D)_{i,j=1\ldots p}$ and set 
$D_k= \tr D^k := \F{\P^r}{\P y^{i_1}_{~i_2} \P y^{i_2}_{~i_3}\ldots\P 
y^{i_k}_{~i_1}}$. Now assume $\phi_\nu \equiv 0$ for $\nu  < n$ thus also 
equations  (\ref{M1M0}). Then we find successively 
\begin{equation}\label{phirek}
 \phi_n = - \F{1}{(-\mu)^k} \sum_{r=k+1}^{2k}\F{1}{r!}\sum_{{\alpha\in 
 Con(S_r)\atop \alpha_1=0~\#\alpha=k}} h_\alpha D_2^{~\alpha_2}\ldots 
 D_r^{~\alpha_r}\GA G_{(n-k)},x-\mu\DA ~.  
\end{equation}
This formula can be proved by induction using the above formulas. Using it  
for $k=n$ shows that $\phi_n$ vanishs identically and the proof is 
completed. \hspace*{\fill}$\square$

A final lemma will be useful in the next section:
\begin{lemma}\label{fFsplit}
Let $F=PF+\GA G\stackrel\stt{,}x-\mu\DA$ be the unique decomposition of the 
invariant function $F$ and let $f$ be homogeneous $Pf=f$. Then $f 
F=f\,PF+\GA fG\stackrel\stt{,}x-\mu \DA$ is the decomposition of the 
product $f F$. 
\end{lemma}
 
\NI {\bf Proof:}
For a homogeneous function $f$ holds $z^A_{~i}\F{\P}{\P z^A_{~j}} f=0\quad 
i,j=1\ldots p$. Hence we have $\GA f G\stackrel\stt{,}x-\mu \DA=f\GA  G 
\stackrel\stt{,}x-\mu\DA$. \hspace*{\fill}$\square$ 

\section{The solution for the projective space $\mathbb{C}P^q$}
The first step towards the construction of a star product on the quotient 
manifold for the projective space as well as for the general Grassmann 
manifolds consists in  "homogenizing" the expression for the Wick product 
of two homogeneous functions. The summands in the series are written as 
products of homogeneous differential operators -- differential operators 
that assign two homogeneous functions another homogeneous function and that 
are thus well defined on the quotient manifold -- and powers of a certain 
nonhomogeneous function. According to the lemma in the last paragraph it 
suffices to calculate the projection for this certain function and its 
powers. In case of the projective space $\mathbb{C}P^q$ the homogeneization 
is simple and was also used in \cite{BBEW95} to find an explicit expression 
for the star product on $\mathbb{C}P^q$. 
\begin{eqnarray}
  \NN F~\stt~ G  & = &  \sum_{r=0}^\infty \F{\lambda^r}{r!} 
  x^{-r}\underbrace{\G(x^r \F{\P^r F}{\P z^{A_1}\ldots\P 
  z^{A_r}}\cdot\F{\P^r G}{\P{z^\Ad}_{A_1}\ldots {z^\Ad}_{A_r}}\D)} \\ & = &  
  \label{homprocCPq}\sum_{r=0}^\infty \F{\lambda^r}{r!} x^{-r} 
  \hspace{12ex}{M_r(F,G)} 
\end{eqnarray} 
Obviously $x\F{\P}{\P z^A}\otimes\F{\P}{\P {z^\Ad}_A}$ is a homogeneous 
differential operator, assigning two homogeneous functions the tensor 
product of homogeneous functions, so its power (followed by multiplication) 
$M_r$ are as well homogeneous. Now a star product for homogenous functions 
$F,G\in C_{hom}^\infty(\mathbb{C}^{q+1}\backslash\{0\})$ is easily 
constructed 
\begin{equation}\label{F*G}
  F * G = P(F~\stt~G)=\sum_{r=0}^\infty P(x^{-r}M_r(F,G))=\sum_{r=0}^\infty 
  P(x^{-r})M_r(F,G)\quad,
\end{equation}
where the last equation sign holds because of lemma \ref{fFsplit}. 
Following  the program described in section \ref{secdec} we obtain 
\[\GA H(x^{-r}),x-\mu\DA  =  x^{-r} -\B{\rho}^* i^* x^{-r}= 
x^{-r}-\mu^{-r}  = -(x-\mu)\sum_{\nu=1}^{r}x^{-\nu}\mu^{\nu-r-1}\]        
In this  case $H(F)$ is a uniquely determined invariant function. 
\[  \Delta(x^{-r})  =  z^A\F{\P}{\P z^A} H(x^{-r}) = x \F{\P}{\P x}H(x^{-r}) = 
  \sum_{\nu=1}^r\nu x^{-\nu} \mu^{\nu-r-1}\]                               
Applying the general recursion formula 
\begin{equation}\label{Rekursion}
  P(F)= \B{\rho}^* i^* F -\lambda P(\Delta F) 
\end{equation}
to the cases $F=x^{-r}$ and $F=\mu^{-1}x^{-r+1}$ as well as to$F=x^{-1}$ 
yields a recursion formula and its beginning:
\[P\G(x^{-r}\D)=\F{1}{\mu +\lambda r} P\G( x^{-(r-1)}\D)\qquad\qquad 
  P\G(x^{-1}\D)=\F{1}{\mu+\lambda} \] This proves the simple \vspace{-1ex} 
\begin{lemma} \label{Px^-r}$~$ 
\vspace{-3ex}\[ P\G(x^{-r}\D)=\prod_{s=1}^r \F{1}{\mu+\lambda s}= 
\lambda^{-r}\F{1}{[c]_r}~,\]
 where \vs{-3ex} 
\[  [y]_r:=y(y+1)\ldots (y+r-1)~,~~ [y]_0:=1~~\mbox{and } c:= 
  \mu\lambda^{-1}+1~.\vspace{-3.7ex}\]  \hspace*{\fill}$\square$ 
\end{lemma}            
In the same way as described here one could calculate $ 
P(x^r)=\prod_{s=0}^{r-1} (\mu-\lambda s)  ~.$  
\begin{corollary}\label{CPq*}
There exists a star product on $G_{1,q}(\mathbb{C})=\mathbb{C}P^q$ defined 
by ($f,g 
\in C^\infty (\mathbb{C}P^q)$) 
\begin{equation}\label{CPn*}
   \pi^* (f*g) :=P(\pi^* f~\stt~ \pi^* g) = \sum_{r=0}^\infty 
   \F{1}{r![c]_r}M_r(\pi^*f,\pi^*g)\quad, 
\end{equation}
where the homogeneous operators $M_r$ are defined by
\begin{equation}\label{CPnM_r}
    M_r(F,G):=x^r \F{\P^r F}{\P z^{A_1}\ldots\P z^{A_r}}\cdot\F{\P^r G}{\P 
    {z^\Ad}_{A_1}\ldots\P {z^\Ad}_{A_r}}=\mu^r\B{\rho}^*i^*\G( \F{\P^r 
    F}{\P z^{A_1}\ldots \P z^{A_r}}\cdot\F{\P^r G}{\P {z^\Ad}_{A_1}\ldots 
    \P {z^\Ad}_{A_r} } \D)\hspace*{\fill}\square
\end{equation}
\end{corollary}

At this point it seems to be only a matter of curiosity that for this star 
product not only the limit $\lambda \rightarrow 0$, but also the opposite 
one, $\lambda \rightarrow \infty$ is well defined as a formal power series.  
And at a glance the space of functions for which this product converges is 
nonempty. This strange possibility of a "strong quantum limit" will be 
preserved in the Grassmann case. 

\section{The projection for Grassmann manifolds} We have to find a 
transfered version of the homogeneization procedure (\ref{homprocCPq}), 
since $x$ is no longer a scalar function. First note that according to the 
usual properties differentiation lowers the degree of  homogeneity by one, 
so for a homogeneous function $F\in\Chom$ which is more precisely a 
homogeneous of degree 0 holds: 
\[
  F(zg)  = F(z)\quad \forall g\in Gl(p,\mathbb{C})\qquad \Longrightarrow 
  \qquad  g^i_{~j}\F{\P F}{\P z^A_{~j}}(zg) = \F{\P F}{\P z^A_{~i}}(z) 
  \qquad { i=1\ldots p \atop A=1\ldots p+q} 
\]
Thus the expression $z^B_{~i}\F{\P F}{\P z^A_{~i}}$ is again homogeneous 
and this suggests to rewrite the expression for the Wick product for the 
case of homogeneous functions $F,G \in \Chom $ in the following way: 
\begin{eqnarray}
 \!\!\NN F~ \stt~ G & = & m\circ \exp 
 \G(\lambda\delta^A_{~B}\delta^j_{~i}\F{\P}  
  {\P z^A_{~i}}\otimes \F{\P}{\P {z^\Ad}^j_{~B}}\D)(F\otimes G)\\   & = & 
  m\circ \exp\G(\lambda S^D_{~C} z^C_{~i}\F{\P}{\P z^A_{~i}}\otimes\F{\P}{\P 
  {z^\Ad}^j_{~A}}{z^\Ad}^j_{~D}\D)(F\otimes G)~=~ 
 \sum_{r=0}^\infty \F{\lambda^r}{r!}\bigg( S^{\otimes r}, 
   M_r(F,G)\bigg) 
\end{eqnarray}
where \vspace{-5ex}
\begin{eqnarray*}
\qquad\qquad S  & = & zx^{-2}z^\Ad\quad, \\ 
M_r(F,G)^{C_1\ldots C_r}{}_{D_1\ldots
D_r} & 
= & z^{C_1}_{~i_1}\ldots z^{C_r}_{~i_r}\F{\P^r F}{\P z^{A_1}_{~i_1}\ldots 
\P 
  z^{A_r}_{~i_r}}\F{\P^r G}{\P {z^\Ad}^{j_1}_{A_1}\ldots 
  {z^\Ad}^{j_r}_{A_r}}{z^\Ad}^{j_1}_{~D_1}\ldots {z^\Ad}^{j_r}_{~D_r}~,\\ 
  m(F\otimes G) & = & F\cdot G~,  
\end{eqnarray*}
and $(\cdot,\cdot)$ is the extension of the usual hermitian product $(A,B)= 
\mbox{tr} A^\Ad B$ to tensor products. The matrix $ S $ is obviously 
hermitian and in contrast to the momentum mapping $J$, which is only 
equivariant with respect to the $U(p)$-action, even $U(p)$-invariant. The 
homogeneous differential operators take values in the $r$-fold tensor 
product of the defining representation of $U(p+q)$. 

Formally the construction of a star product for homogeneous functions on 
$\mathbb{C}^{(p+q)\cdot p~*}$ is completely analogous to the 
$\mathbb{C}P^q$-case of the last section
\begin{equation}\label{F*G Gras}
  F*G:= P(F~\stt ~G)= \sum_{r=0}^\infty \F{\lambda^r}{r!} P\bigg( 
  S^{\otimes r},M_r(F,G)\bigg) = \sum_{r=0}^\infty \F{\lambda^r}{r!}\bigg( 
  P(S^{\otimes r}), M_r(F,G)\bigg ) ~, 
\end{equation}
where again the last equation sign is justified by lemma (\ref{}). Yet the 
evaluation of the projection $P(S^\otr)$ shows a new aspect for which we 
introduce the notion of the action of the symmetric group on the tensor 
product space$(\mathbb{C}^s)^\otr$:
\begin{definition}\label{defrho}
For $\sigma \in S_r$ let $\rho(\sigma):{(\mathbb{C}^s)}^\otr 
\longrightarrow {(\mathbb{C}^s)}^\otr$ be the representation defined by 
linear extension of the prescription $ 
  \rho(\sigma)(e_{A_1}\otimes\ldots\otimes\  e_{A_r}) = e_{A_{\sigma(1)}} 
  \otimes \ldots \otimes e_{A_{\sigma(r)}}\qquad 1\leq A_i\leq s$, where 
$\{e_A\}_{A=1\ldots s}$ is a basis of $\mathbb{C}^s$. This representation 
can linearly be extended to the group algebra. As well let for $S \in 
Gl(s,\mathbb{C})$ be $D(S):{(\mathbb{C}^s)}^\otr \longrightarrow 
{(\mathbb{C}^s)}^\otr$ the representation defined by $ 
  D(S):=S^\otr$, which in a basis reads $D(S)(e_{A_1}\otimes\ldots\otimes\ 
  e_{A_r}):= e_{B_1}\otimes\ldots\otimes e_{B_r} S^{B_1}_{~A_1} \ldots 
  S^{B_r}_{~A_r}.$  
\end{definition}
Obviously the two representations commute, that means $\rho(\sigma)\circ 
D(S)=D(S)\circ \rho(\sigma)$ for all $ 
\sigma\in S_r$ and $S\in Gl(s,\mathbb{C})$.
Linear extension of $D$ to a representation of the endomorphism algebra 
$End(\mathbb{C}^s)$ leads to the subalgebra of all mappings that commute 
with permutations. We will denote this subalgebra as 
$\mathbb{C}[D(Gl(s,\mathbb{C}))]$.
 
Using the notation introduced in lemma \ref{Zerlegung} concerning the 
conjugacy classes of the symmetric group we can formulate the essential 
\begin{lemma}\label{PS^r} For the 
projection $P(S^\otr)$ the following recursion relation holds:
\begin{equation}\label{RekPS^r}
  \lambda^r \rho\Big(\sum_{\alpha\in Con(S_r)} c^{|\alpha|} k_\alpha \Big) 
  P(S^\otr)= T^\otr~, 
\end{equation}
where $T:=\mu\B{\rho}^* i^* S=\mu^{-1}\zeta^\Ad \zeta,~ 
c:=\lambda^{-1}\mu+p$ and $k_\alpha:=\sum_{\sigma\in\alpha} \sigma$ is 
defined to be the sum over all permutations of the conjugation class 
$\alpha$ in the group algebra $\mathbb{C}[S_r]$. 
\end{lemma}

\NI {\bf Proof:} The formula will be proved by induction. 
For any constant hermitian $p+q$ matrix $\Phi$ holds:
\[
  \G\GA H\Big((S,\Phi)\Big),x-\mu\D\DA  =  
  -\mu^{-1}(zx^{-1}(x-\mu)x^{-1}z^\Ad,\Phi) 
   = -\mu^{-1}\GA x^{-1}z^\Ad\Phi zx^{-1},x-\mu\DA 
\]
and we can choose 
\begin{equation}\label{H(S,Phi)}
  H\Big((S,\Phi)\Big)=-\mu^{-1}z^\Ad\Phi z x^{-1}~.
\end{equation}
The calculation of $\Delta(S,\Phi)$ is simplified by the fact that 
$z^A_{~i}\F{\P}{\P z^A_{~j}}$ is the fundamental vector field for the 
holomporphic $Gl(p,\mathbb{C})$-action, i.e. $((z,z^\Ad),g)\longrightarrow 
(zg,z^\Ad)$ under which the term $zx^{-1}$ is invariant. 
\[ \Delta\Big((S,\Phi)\Big)=z^A_{~i}\F{\P}{\P z^A_{~j}} H^i_{~j}(S,\Phi)=
\mu^{-1} p 
(S,\Phi)\qquad\Longrightarrow\qquad   \Delta(S)=\mu^{-1}p S\] Using the 
recursion formula (\ref{Rekursion}) one finds $\lambda c P(S) 
= T$, which proves the start of the induction. Now assume the formula is 
proved for $i<r$. A possible choice for $H$ is 
\[ H(S^{\otimes r}) = S^{\otimes r-1}\otimes H(S) + 
H(S^{\otimes r-1})\otimes\mu^{-1}T \] and it follows  
\[ 
\begin{array}{rcccccl}
 \DS \Delta(S^\otr) & \DS = &\DS S^{\otimes r-1}\otimes\Delta(S)&\DS + 
 &\DS(z^A_{~i}\F{\P}{\P 
  z^A_{~j}} S^{\otimes r-1} )\otimes H^i_{~j}(S)&\DS  + &\DS \Delta 
  (S^{\otimes r-1}) \otimes \mu^{-1} T 
\\   &\DS = &\DS \mu^{-1}p S^{\otimes r} &\DS + &\DS \sum_{s=1}^{r-1} 
\rho(\tau_{sr})S^\otr &\DS + &\DS \mu^{-1} \Delta(S^{\otimes r-1})\otimes T~,
\end{array}\]
where $\tau_{sr}$ is the transposition of $s$ and $r$. Comparing the 
recursion formula (\ref{Rekursion}) analogously to the proof of lemma 
(\ref{Px^-r}) for $P(S^\otr)$ and $\mu^{-1}P(S^{\otimes r-1})\otimes T$ 
taking into consideration that $P$ and $\rho$ commute, yields: 
\[ \G(1+\lambda\mu^{-1} p + \lambda\mu^{-1}\sum_{s=1}^{r-1}\rho(\tau_{ir})\D) 
  P(S^\otr) = \mu^{-1} P(S^{\otimes r-1})\otimes T \]                       
Using the induction hypothesis and the relation (\ref{step}) leads to
\[ \lambda^r\bigg(c + \sum_{s=1}^{r-1} \rho (\tau_{sr})\bigg)
\Bigg(\sum_{\sigma\in 
S_{r-1}\subset S_r \atop \sigma(r)=r} 
c^{r-1-\#[\sigma]}\rho(\sigma)\Bigg)P(S^\otr)=\lambda^r \sum_{\sigma\in 
S_r}c^{r-\#[\sigma]}\rho(\sigma)P(S^\otr)= T^\otr ~, 
\] 
which proves the lemma.\hspace*{\fill}$\square$ 

\NI {\bf Remark:} In the same way one can show that for the positive powers 
holds
\begin{equation}\label{PzAdz^r}
  P\G({z z^\Ad}^\otr\D) = \rho\Bigg(\sum_{\alpha \in 
  Con(S_r)}(-\lambda\mu^{-1})^{\#\alpha} k_\alpha\Bigg) 
  \G({\zeta\zeta^\Ad}^\otr\D)~.
\end{equation}
 
Since $\{k_\alpha\}_{\alpha\in S_r}$ is a basis of the center of the group 
algebra $\mathbb{C}[S_r],~ \sum_{\alpha\in Con(S_r)} c^{|\alpha|}k_\alpha$ 
is an element of the center which will be invertible for generic $c$. For 
the inverse element we can write 
\begin{equation}\label{s}
  \bigg(\sum_{\alpha} c^{|\alpha|}k_{\alpha}\bigg 
  )^{-1}=:\sum_{\alpha}s_{\alpha}(c)k_{\alpha}~, 
\end{equation}
where $s_\alpha$ are rational functions which are to be determined. First 
we formulate a trivial corallary of lemma \ref{PS^r}. 
\begin{corollary}\label{corgras*}
There exists a star product on $G_{p,q}(\mathbb{C})$ defined by ($f,g \in 
C^\infty(G_{p,q}(\mathbb{C}))[[\lambda]]$) 
\begin{equation}\label{gras*}
\pi^*(f\:*\:g)  = P(\pi^* f~\stt ~\pi^*g ) = \sum_{r=0}^\infty  
\F{\mu^r}{r!}\B{\rho}^*i^*\G\GA \rho\bigg(\sum_{\alpha\in Con(S_r)} 
s_\alpha(c) k_\alpha\bigg)
\F{\P^r\pi^* f}{\P z^{A_1}\ldots \P z^{A_r}},\F{\P^r\pi^* g}{\P 
{z^\Ad}_{A_1}\ldots {z^\Ad}_{A_r}}\D\DA, 
\end{equation}
where the rational functions $s_\alpha$ are defined by (\ref{s}) and 
$c=\lambda^{-1}\mu+p$. Here we regard the derivatives $\F{\P^r F}{\P 
z^{A_1}_{~i_1}\ldots\P z^{A_r}_{~i_r}}$ as $(i_1,\ldots,i_r)$-component of 
a tensor in $(\mathbb{C}^p)^\otr$ and $\GA\cdot,\cdot\DA$ denotes the 
hermitian inner product of $(\mathbb{C}^p)^\otr$. 
\hspace*{\fill}$\square$ 
\end{corollary} 
A comparison with corollary \ref{CPq*} yields as a condition for 
$s_\alpha$: 
\begin{equation}\label{sums_alpha}
  \sum_{\sigma\in S_r} s_{[\sigma]}(c)=\sum_{\alpha\in Con(S_r)}h_\alpha 
  s_\alpha(c) = \F{1}{[c]_r}~. 
\end{equation}

\section{The irreducible differential operators and their coefficients}
This section is devoted to the determination of the rational functions 
$s_\alpha$ defined by (\ref{s}) that arise in the star product on Grassmann  
manifolds (\ref{gras*}). According to Wedderburn's theorem the group 
algebra $\mathbb{C}[G]$ of a finite group $G$ decomposes directly into 
simple endomorphism rings $End(\mathbb{C}^{n_a}),~a=1\ldots f$ of dimension 
$n_a^2$ such that $\sum n_a^2$ equals the order $|G|$ of the group. Each 
simple ring $End(\mathbb{C}^{n_a})$ may be seen as an isotypical submodule 
of the regular representation of the group $G$ on its group algebra 
$\mathbb{C}[G]$ obtained by left multiplication. These isotypical 
submodules $End(\mathbb{C}^{n_a})$ are  associated to a certain irreducible 
representation $\rho^a$, which is therein contained with multiplicity 
$n_a$. To sum up, all irreducible representations  $\rho^a$ of a finite 
group $G$ are contained in the group algebra and their multiplicity equals 
their dimension $n_a$. The number $f$ of inequivalent irreducible 
representations equals the number of conjugacy classes, which is in case of 
the symmetric group $S_r$ the number of partitions of $r$. So the center of 
the group ring $\mathbb{C}[G]$ has two natural bases, the conjugation 
classes $k_\alpha$ and the units $e_{a}$ of the simple rings 
$End(\mathbb{C}^{n_a})$ into which the group algebra decomposes. Between 
both there is a simple base transformation in terms of characters 
\cite{Boe67}, which for the case of the symmetric group $S_r$ reads 
\begin{equation}\label{k-e-basistrafo}
  e_a=\F{n_a}{r!}\sum_{\alpha=1}^f\chi_\alpha^ak_\alpha\qquad\qquad 
  k_\alpha=h_\alpha\sum_{a=1}^f \F{\chi_\alpha^a}{n_a} e_a 
\end{equation}
Here $\chi_\alpha^a$ is the character of the irreducible representation 
$\rho^a$ evaluated at an element of the conjugacy class $\alpha$. Note that 
for the symmetric group a group element and its inverse are in the same 
class, such that the character is real. Of course the basis 
$\{e_a\}_{a=1\ldots f}$ is more appropriate for considering the inverse 
(\ref{s}): 
\begin{equation}\label{tdef}
  \bigg(\sum_{\alpha} c^{|\alpha|}k_\alpha\bigg)^{-1}= 
  \sum_{a}\F{n_a}{\sum_{\alpha} h_\alpha\chi_\alpha^a 
  c^{|\alpha|}}e_a=:\sum_a \F{1}{t_a(c)} e_a\qquad 
  t_a(c):=\F{1}{n_a}\sum_\alpha h_\alpha\chi_\alpha^a c^{|\alpha|} 
\end{equation}
The relation between the rational functions $s_\alpha$ introduced in 
(\ref{s}) and the polynomials $t_a$ is 
\begin{equation}\label{s-t}
  s_\alpha(c)=\sum_a\F{n_a\chi_\alpha^a}{r!}\F{1}{t_a(c)} 
\end{equation}
As noted above any isotypical submodule $End(C^{n_a})$ of the group ring or 
equivalently any irreducible representation $\rho^a$ (up to isomorphism)   
or any unit $e_a$ is characterised by a partition, which is usually 
represented by a frame. We quote the dimension formula \cite{Boe67}: 
\begin{theorem}
The dimension $n_{[m]}$ of an irreducible representation $\rho^{[m]}$ of 
$S_r$ characterised by a frame $[m]=[m_1,\ldots m_k]$ with $k$ rows of  
lengths $m_1\geq m_2\geq\ldots m_k\geq 0,~|m|:=\sum_{i=1}^k m_i=r$, and 
thus the multiplicity in the group algebra is given by 
\begin{equation}\label{n_a}
  n_{[m]}=r!\F{\prod_{i<j} (l_i-l_j)}{l_1!\ldots l_k!}\qquad l_i=m_i+k-i 
\end{equation}
Note that the dimension is of course invariant if one adds certain rows of 
length zero, so the right hand side is in fact unchanged under 
$[m]\longmapsto[m^\prime]=[m,0]$.                                              
\end{theorem}
Now we consider the decomposition of the $r$-fold product of the vector 
space $\mathbb{C}^s$ into symmetry classes. $(\mathbb{C}^s)^\otr$ carries 
the  representation $\rho$ defined above of the symmetric group $S_r$ and 
is thus decomposable into irreducible components characterised by a frame. 
The  irreducible components may be contained with multiplicity $d_a$ 
(possibly 0). 
\begin{equation}\label{C^sr,S_r}
  (\mathbb{C}^s)^\otr=d_1\mathbb{C}^{n_1}\oplus\ldots\oplus 
  d_f\mathbb{C}^{n_f}\qquad \mbox{w.r.t. }S_r   
\end{equation}
Since the action $D$ of $Gl(s,\mathbb{C})$ does commute with the action of 
the symmetric group Schur's lemma states, that an element of 
$Gl(s,\mathbb{C})$ can only act as a multiple of unity between the 
irreducible summands of (\ref{C^sr,S_r}) which has to be zero between 
$\mathbb{C}^{n_i}$ and $\mathbb{C}^{n_j}$ if $i\neq j$. This leads to 
another decomposition 
\begin{equation}\label{C^sr,Gl(C)}
  (\mathbb{C}^s)^\otr = n_1 \mathbb{C}^{d_1}\oplus\ldots\oplus 
  n_f\mathbb{C}^{d_f} \qquad \mbox{w.r.t. }Gl(s,\mathbb{C})~, 
\end{equation}
where here $\mathbb{C}^{d_a}$ are irreducible representation spaces  of an 
irreducible representation $D^a$ of $Gl(s,\mathbb{C})$, which arise with 
multiplicity $n_a$. As a trivial consequence the isotypical components are 
identical for both actions:
\begin{equation}\label{C^sr,iso}
  (\mathbb{C}^s)^\otr = \mathbb{C}^{d_1}\otimes\mathbb{C}^{n_1}\oplus\ldots
  \oplus \mathbb{C}^{d_f}\otimes \mathbb{C}^{n_f}
\end{equation}
The units $e_a$ of the simple components of the group ring 
$\mathbb{C}[S_r]$ act as projectors into the isotypical components here as 
well. A unit $e_a$ associated to a frame whose number of rows exceeds $s$, 
acts as an annihilator, $\rho(e_a)=0$, since the degree of 
antisymmetrisation is too large, so $n_a$ must be 0 in this case. So the 
representations $D^a,~a=1\ldots r$ of $Gl(s,\mathbb{C})$ are as well 
characterised by frames. We quote again the dimension formula \cite{Boe67}: 
\begin{theorem}
For a frame with more than $s$ rows of nonzero length the dimension of the 
associated representation $D^{[m]}$ of $Gl(s,\mathbb{C})$ is zero. So let 
$[m]$ be a frame with $s$ rows of length $m_1\geq \ldots\geq m_s\geq 0$, 
where the   lowest rows may be of length zero. Then the dimension $d_{[m]}$ 
of the associated representation $D^a$ of $Gl(s,\mathbb{C})$ is 
\begin{equation}\label{d_a}
  d_{[m]}=\F{\prod_{i<j}(l_i - l_j)}{(s-1)!(s-2)!\ldots 1!0!}\qquad 
  l_i=m_i+s-i~.
\end{equation}
\end{theorem}

Finally we recapitulate the considerations leading to Frobenius' formula 
relating the characters of the respective representations $\rho^a$ and 
$D^a$ of $S_r$ and $Gl(s,\mathbb{C})$. Consider the composite actions of 
$A\in Gl(s,\mathbb{C})$ and $\sigma\in S_r$  on $(\mathbb{C}^s)^\otr$ and 
assume $\sigma$ belongs to the conjugacy class $\alpha$. Using the 
isotypical decomposition (\ref{C^sr,iso}) the mapping decomposes as 
\[ D(A)\circ\rho(\sigma)=\rho(\sigma)\circ D(A)=D^1(A)\otimes \rho^1(\sigma)
\oplus \ldots \oplus D^f(A)\otimes \rho^f(\sigma) \] 
and for its trace one finds on the one hand 
\[ \tr\Big(D(A)\circ\rho(\sigma)\Big)=\sum_{a=1}^f \tr D^a(A)\cdot\tr
 \rho^a(\sigma)= \sum_{a=1}^f \phi(A)\chi_\alpha^a~, \]                    
where $\phi^a$ is  the character of the representation $D^a$, and on the 
other hand according to definition \ref{defrho} 
\[\tr\Big(D(A)\circ\rho(\sigma)\Big)=A^{i_1}_{~i_{\sigma(1)}}\cdot\ldots\cdot
A^{i_r}_{~i_{\sigma(r)}}=a_1^{~\alpha_1}\cdot\ldots\cdot 
a_r^{~\alpha_r}~,\] where $a_k:=\tr A^k$ is the trace of the powers of $A$. 
Finally using the orthogonality relation, which may be read from equations 
(\ref{k-e-basistrafo}) yields Frobenius' 
\begin{theorem}
The relation between the characters of the irreducible representations 
$D^{[m]}$ of $Gl(s,\mathbb{C})$ and $\rho^{[m]}$ of $S_r$ associated to the 
same frame $[m]=[m_1,\ldots,m_s],~|m| =r$ is given by
\begin{equation}\label{Frobformel}
  \phi^{[m]}(A)=\F{1}{r!}\sum_{\alpha\in Con(S_r)} h_\alpha 
  a_1^{~\alpha_1}\cdot\ldots\cdot a_r^{~\alpha_r}  \chi_\alpha^{[m]}~,
\end{equation}
where $a_i:=\tr A^i$ and $h_\alpha$ is the number of element in the 
conjugacy class $\alpha$.\hspace*{\fill}$\square$
\end{theorem}
Summing up all informations the proof of the following corollary is 
trivial. 
\begin{corollary}
The polynomial $t_{[m]}$ introduced by (\ref{tdef}) associated to 
the frame $[m]$ is given by 
\begin{equation}\label{tsol}
  t_{[m]}(c)=[c]_{m_1}[c-1]_{m_2}\cdot\ldots\cdot [c-p+1]_{m_p}
\end{equation}
\end{corollary}
{\bf Proof:} It suffices to know the values of the polynomial for all 
natural values $s\geq p$. Then 
\[ t_{[m]}(s)=\F{r!}{n_{[m]}}\phi^{[m]}({\bf 1}_s)=[s]_{m_1}[s-1]_{m_2}\cdot
\ldots\cdot [1]_{m_s}
 \]                                                                          
But the frames under consideration in our case have at most $p$ rows, since 
otherwise $e_{[m]}$ acts as annihilator on tensors in 
$(\mathbb{C}^p)^\otr$. Thus in our case $m_{p+1}=\ldots m_s=0$ and the 
formula is proved.\hspace*{\fill}$\square$

This formula has a simple graphical interpretation: Consider an arbitrary 
frame:

{\tiny\hspace{-10ex} 
\unitlength 0.80mm
\linethickness{0.4pt}
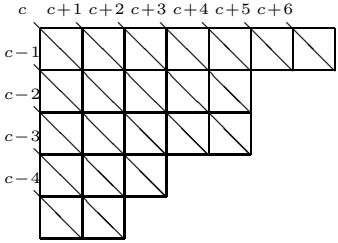
\begin{picture}(110.00,50.00)(00.00,60.00)
\put(10.00,100.00){\line(1,0){49.00}}
\put(59.00,100.00){\line(0,-1){7.00}}
\put(59.00,93.00){\line(-1,0){49.00}}
\put(10.00,86.00){\line(1,0){35.00}}
\put(10.00,79.00){\line(1,0){35.00}}
\put(45.00,79.00){\line(0,1){21.00}}
\put(10.00,72.00){\line(1,0){21.00}}
\put(31.00,72.00){\line(0,1){28.00}}
\put(10.00,65.00){\line(1,0){14.00}}
\put(24.00,65.00){\line(0,1){35.00}}
\put(10.00,65.00){\line(0,1){35.00}}
\put(17.00,65.00){\line(0,1){35.00}}
\put(38.00,79.00){\line(0,1){21.00}}
\put(52.00,93.00){\line(0,1){7.00}}
\put(7.00,103.00){\makebox(0,0)[cc]{$c$}}
\put(21.00,103.00){\makebox(0,0)[cc]{$c\!+\!2$}}
\put(7.00,96.00){\makebox(0,0)[cc]{$c\!-\!1$}}
\put(7.00,89.00){\makebox(0,0)[cc]{$c\!-\!2$}}
\put(14.00,103.00){\makebox(0,0)[cc]{$c\!+\!1$}}
\put(28.00,103.00){\makebox(0,0)[cc]{$c\!+\!3$}}
\put(35.00,103.00){\makebox(0,0)[cc]{$c\!+\!4$}}
\put(42.00,103.00){\makebox(0,0)[cc]{$c\!+\!5$}}
\put(49.00,103.00){\makebox(0,0)[cc]{$c\!+\!6$}}
\put(7.00,82.00){\makebox(0,0)[cc]{$c\!-\!3$}}
\put(7.00,75.00){\makebox(0,0)[cc]{$c\!-\!4$}}
\put(9.00,73.00){\line(1,-1){8.00}}
\put(9.00,80.00){\line(1,-1){15.00}}
\put(9.00,87.00){\line(1,-1){15.00}}
\put(9.00,94.00){\line(1,-1){22.00}}
\put(9.00,101.00){\line(1,-1){22.00}}
\put(16.00,101.00){\line(1,-1){22.00}}
\put(23.00,101.00){\line(1,-1){22.00}}
\put(30.00,101.00){\line(1,-1){15.00}}
\put(37.00,101.00){\line(1,-1){8.00}}
\put(51.00,101.00){\line(1,-1){8.00}}
\put(44.00,101.00){\line(1,-1){8.00}}
\put(135.00,86.00){\makebox(0,0)[cc]{{\large $\longmapsto$}{\small $
\qquad t_{[m]}(c)=(c+6)(c+5)(c+4)(c+3)^2(c+2)^3\times$}}}
\put(155.00,79.00){\makebox(0,0)[cc]{{\small $(c+1)^3c^3(c-1)^3(c-2)^2
(c-3)^2(c-4)$}}}
\end{picture}}

\NI For any box lying on the main diagonal $t_{[m]}$ contains a factor $c$,
 for
any box on the first upper resp. lower diagonal it contains a linear factor  
$c+1$ resp $c-1$ and so forth. The relation  (\ref{sums_alpha}) is now 
recognized as a direct consequence of equation (\ref{s-t}) and the 
orthogonality relation. The star product on the Grasmann manifolds 
$G_{p,q}(\mathbb{C})$ of corollary \ref{corgras*} then takes the form:  
\begin{theorem}\label{G_pq*}
\hspace{-1.5ex} There exists a star product on $G_{p,q}(\mathbb{C})$ such that
 for $f,g\! 
\in \!C^\infty (G_{p,q}(\mathbb{C}))[[\lambda]]$: 
\begin{eqnarray}
\NN  \!\!\pi^*(f*g) & = & P(\pi^* f~\stt~\pi^*g)\\
\label{finales*}  & = & \sum_{r=0}^\infty \F{\mu^r}{r!}\sum_{|m|=r \atop 
  m_{p+1}=0}\F{1}{t_{[m]}(c)}\B{\rho}^*i^*\G\GA \rho(e_{[m]})\F{\P^r \pi^* 
  f}{\P z^{A_1}\ldots \P z^{A_r}},\rho(e_{[m]})\F{\P^r\pi^* g}{\P 
  {z^\Ad}_{A_1}\ldots \P {z^\Ad}_{A_r}}\D\DA 
\end{eqnarray}
Here $\rho(e_{[m]})$ denotes the projector of $(\mathbb{C}^p)^\otr$ onto 
the symmetry class characterised by the frame $[m]$ and $t_{[m]}(c)
 =[c]_{m_1}[c-1]_{m_2}\cdot\ldots\cdot[c-p+1]_{m_p}$. In the sum over the 
frames only such frames yield a nonzero contribution whose number of rows 
is not greater than $p$. \hspace*{\fill}$\square$
\end{theorem}
Obviously this theorem contains Corollary (\ref{CPq*}) as a special case 
for $p=1$. The product remains well defined in the limit $\lambda 
\rightarrow \infty$, which is here even more remarkable, since frames 
containing one factor $(c+p-1)$ do arise, but those containing the factor 
$c-p=\mu\lambda^{-1}$ do not. Karabegov already stated in Theorem 3 
\cite{Kar94} that for representative functions $f,g$ the dependence of 
$f*g$ on $\lambda^{-1}$ is rational with no poles for $\lambda\rightarrow 
0$. For arbitrary smooth functions formula (\ref{finales*}) states that the 
dependence on $\lambda$ of the product is given by an infinite sum of 
rational functions, obviously in accordance with Karabegov's result. By a 
desription of the star product via representative functions as it is done 
for the projective space in \cite{BBEW96} one could verify Karabegov's 
proposition directly for   Grassmann manifolds. The expected isomorphy of 
the star algebras $(C^\infty(G_{p,q}(\mathbb{C})),*)$ and 
$(C^\infty(G_{q,p}(\mathbb{C})),*)$  can probably as well be established 
most easily in terms of reprsentative functions. This would extend the 
duality of $G_{p,q}(\mathbb{C})$ and $G_{q,p}(\mathbb{C})$ to their star 
algebras.  

%
%

\NI {\bf Acknowledgement:} I would like to thank M. Bordemann, whose idea of 
the projection method was actually the starting point of this work.

\end{document}